\documentclass[12pt,dvips]{amsart}
\usepackage{amsfonts, amssymb, latexsym, epsfig, amscd}

\newtheorem{Definition-Proposition}[Theorem]{Definition-Proposition}
\newtheorem{Main Conjecture}[Theorem]{Main Conjecture}

\theoremstyle{remark}

\theoremstyle{plain}

\newcommand{\cellsize}{13}
\newlength{\cellsz} 
\newsavebox{\cell}
\sbox{\cell}{\begin{picture}(\cellsize,\cellsize)
\put(0,0){\line(1,0){\cellsize}}
\put(0,0){\line(0,1){\cellsize}}
\put(\cellsize,0){\line(0,1){\cellsize}}
\put(0,\cellsize){\line(1,0){\cellsize}}
\end{picture}}
\newcommand\cellify[1]{\def\thearg{#1}\def\nothing{}%
\ifx\thearg\nothing
\vrule width0pt height\cellsz depth0pt\else
\hbox to 0pt{\usebox{\cell} \hss}\fi%
\vbox to \cellsz{
\vss
\hbox to \cellsz{\hss$#1$\hss}
\vss}}
\newcommand\tableau[1]{\vtop{\let\\\cr
\baselineskip -16000pt \lineskiplimit 16000pt \lineskip 0pt
\ialign{&\cellify{##}\cr#1\crcr}}}
\newcommand{\kellsize}{24}
\newlength{\kellsz} 
\newsavebox{\kell}
\sbox{\kell}{\begin{picture}(\kellsize,\kellsize)
\put(0,0){\line(1,0){\kellsize}}
\put(0,0){\line(0,1){\kellsize}}
\put(\kellsize,0){\line(0,1){\kellsize}}
\put(0,\kellsize){\line(1,0){\kellsize}}
\end{picture}}
\newcommand\kellify[1]{\def\thearg{#1}\def\nothing{}%
\ifx\thearg\nothing
\vrule width0pt height\kellsz depth0pt\else
\hbox to 0pt{\usebox{\kell} \hss}\fi%
\vbox to \kellsz{
\vss
\hbox to \kellsz{\hss$#1$\hss}
\vss}}
\newcommand\ktableau[1]{\vtop{\let\\\cr
\baselineskip -16000pt \lineskiplimit 16000pt \lineskip 0pt
\ialign{&\kellify{##}\cr#1\crcr}}}

\newcommand{\sellsize}{36}
\newlength{\sellsz} 
\newsavebox{\sell}
\sbox{\sell}{\begin{picture}(\sellsize,20)
\put(0,0){\line(1,0){\sellsize}}
\put(0,0){\line(0,1){\sellsize}}
\put(\sellsize,0){\line(0,1){\sellsize}}
\put(0,\sellsize){\line(1,0){\sellsize}}
\end{picture}}
\newcommand\sellify[1]{\def\thearg{#1}\def\nothing{}%
\ifx\thearg\nothing
\vrule width0pt height\sellsz depth0pt\else
\hbox to 0pt{\usebox{\sell} \hss}\fi%
\vbox to \sellsz{
\vss
\hbox to \sellsz{\hss$#1$\hss}
\vss}}
\newcommand\stableau[1]{\vtop{\let\\\cr
\baselineskip -16000pt \lineskiplimit 16000pt \lineskip 0pt
\ialign{&\sellify{##}\cr#1\crcr}}}

\begin{document}
\title{The Joseph Greenberg problem: combinatorics\\ and comparative linguistics}

\author{Alexander Yong}
\address{Dept. of Mathematics, University of Illinois at
Urbana-Champaign, Urbana, IL 61801}

\email{ayong@uiuc.edu}

\date{September 30, 2013}
\maketitle


\pagestyle{plain}

\section{Introduction}
In 1957, the eminent linguist Joseph H. Greenberg (1915--2001) proposed the method of \emph{Mass comparison} (also known as \emph{multilateral comparison}) for
determining genetic relatedness between languages \cite{Greenberg1}. I first learned about his work though the
PBS/BBC documentary ``In search of the first language''; a transcript (quoted below) is found
at \textsf{http://www.pbs.org/wgbh/nova/transcripts/2120glang.html}. 

Allowing Greenberg himself to summarize his idea:
\begin{quote}
``I usually had preliminary notebooks in which I took those elements of a language, which, on the whole, we know are the most stable over time. These are things like the personal pronouns, particularly first and second person, names for the parts of the human body... I would look at a very large number of languages in regard to these matters, and I did find that they fell into quite obvious groupings.''
\end{quote}
This ``presorting'' technique is controversial, see, e.g., \cite{Ringe,Greenberg:defense}. Its
use \cite{Greenberg:americas} to see that $650$ languages of North and 
South America fall
into three families (Eskimo-Aleut, Na-Dene and Amerind) is hotly debated. For example, James Matisoff argued:
\begin{quote}
``Eyeballing data is prescientific, or nonscientific. There are so many ways you can be led astray, because very often, words look as if they have some connection, and they have no historical connection whatsoever.''
\end{quote}

Now, a succinct and combinatorial argument for Mass comparison is in \cite{Wiki}:
\begin{quote}
``He also criticized the prevalent view that comprehensive comparisons of two languages at a time (which commonly take years to carry out) could establish language families of any size. He pointed out that, even for $8$ languages, there are already $4,140$ ways to classify them into distinct families, while for $25$ languages there are $4,749,027,089,305,918,018$ ways.'' 
\end{quote}

The two numbers offered are from \cite{Greenberg2} (he gives no other such enumerations). Indeed,
the exactness of his enumeration for $25$ languages makes a strong
visual impression.  However, in the interest of accuracy, we point out that
in fact it is \emph{slightly} erroneous: the correct enumeration is $4,638,590,332,229,999,353$. 


	Combinatorialists will recognize that the numbers Greenberg wanted to be \emph{Bell numbers} $B(n)$. Nowadays, 
a quick lookup at the On-Line Encyclopedia of Integer sequences 
(\textsf{http://oeis.org}) detects the discrepancy.
However, since his count may be of some historical interest, we elaborate upon the correction, and
how one comes to notice it.

Actually, a purpose for this elaboration is pedagogical. The author introduces the ``Joseph Greenberg problem'' to introductory classes in combinatorics: \emph{Compute Greenberg's stated numbers.}
Thus, the discussion might be of interest to the combinatorics instructor, or to the reader who is not already versed on the topic. 

We now turn to the main goal of this paper. While Greenberg used Bell numbers to support Mass comparison,
the same combinatorics does not support the plausibility of the resulting Americas classification. 
The probability, under the uniform distribution, of a classification of $650$ languages having less than $100$ families is near zero, but almost $100\%$ for those in the range $[120,150]$. 

This estimate is close to the more generally accepted viewpoint that there are between $150$ and $180$ families (a list of families is found in, e.g., \cite{Wikib}). In fact, most random classifications with this number of families and languages have a moderate number ($9$ to $19$) of \emph{language isolates}. This is somewhat consistent with the number of isolates/unclassified languages in the actual consensus classification 
(we seem to underestimate the number of isolates/unclassified 
languages by a factor of $5$ to $10$).


We propose these combinatorial tests as simple mathematical baselines to quantitatively measure the typicality of an actual language classification, in cases of
(subjective) uncertainty. Since these tests use a different perspective than more standard statistical methods,
we contribute the analysis to the above debate within comparative linguistics.

\section{Multiset counting, Stirling and Bell numbers, and generating series}

\subsection{The $8$ languages case} Greenberg's calculation of $4,140$ for the number of ways to classify $8$ languages into families is correct. One can first try to determine the number using multiset counting. 

Suppose that the (Native American) languages to be classified are
\[\mbox{Alutiiq, Eyak, Cup'ik, Naukan, Mahican, Inupiaq, Tlinqit, Kalallisut}\]
There are $22$ partitions of $8$.
Each partition describes language family sizes. For example,
in the case $5+2+1$, classifications correspond to to arrangements of the multiset 
 $A,A,A,A,A,B,B,C$. Hence, the arrangement
$A,B,A,A,C,A,B,A$ encodes 
\[A\leftrightarrow\{\mbox{Alutiiq, Cup'ik, Naukan, Inupiaq, Kalallisut}\},  \ B\leftrightarrow\{\mbox{Eyak, Tlinqit}\}, \  C\leftrightarrow\{\mbox{Mahican}\}\]
Thus there are ${8\choose 5 \ 2 \ 1}=\frac{8!}{5!2!1!}$ such arrangements. 

For the partition $4+2+2$, we rearrange $A,A,A,A,B,B,C,C$ and both
\[\mbox{$A,B,B,A,C,C,A,A$ \ \ and \ \ $A,C,C,A,B,B,A,A$}\]
encode the same language groupings. As there are two repeated parts one corrects for this
overcount by dividing by $\frac{1}{2!}$ so that the number of distinct groupings is
$\frac{1}{2!}{8\choose 4 \ 2 \ 2}$. 

Similarly, for the partition $2+2+2+1+1$ the count is $\frac{1}{3!}\frac{1}{2!}{8 \choose 2 \ 2 \ 2 \ 1 \ 1}$. 
It is not too laborious to compute all $22$ numbers of this kind, add them all up, and thus recover Greenberg's stated
number. Perhaps Greenberg did somesuch calculation as a check.

\subsection{The $25$ languages case} 
The method just used for the $8$ language case becomes unpalatable for $25$ languages because
now there are $1,958$ partitions. It is logical to discuss now standard combinatorics, found in textbooks such as \cite{Brualdi}.

Suppose $\{a_n\}_{n\in {\mathbb Z}_{\geq 0}}$ is a sequence of
answers to a counting problem
(such as the Joseph Greenberg problem for $n$ languages). 
The uninitiated might find it suprising that one can ever get traction on 
a problem by considering the infinitely many
such problems, and then rephrasing the question as the
coefficient of $\frac{x^n}{n!}$ in the {\bf exponential generating series}
\[\sum_{n=0}^{\infty} a_n \frac{x^n}{n!}
=a_0 + a_1 \frac{x^1}{1!} + a_2 \frac{x^2}{2!}+a_2\frac{x^3}{3!}+\cdots
.\]


Let $S(n,k)$ be the number of ways to split $n$ languages into exactly $k$ language families; this is known as the {\bf Stirling number}. One has the exponential generating series identity
\begin{equation}
\label{eqn:expgen}
\sum_{k=0}^{\infty}k!S(n,k)\frac{x^n}{n!}=(e^x-1)^k.
\end{equation}
This can be obtained by the {\bf product rule} for exponential generating series: Given 
\[f(x)=\sum_{k=0}^{\infty} f_k \frac{x^k}{k!}, \mbox{ \ and \ }
g(x)=\sum_{\ell=0}^{\infty} g_\ell \frac{x^{\ell}}{\ell!},\]
the coefficient of $\frac{x^n}{n!}$ in $h(x)=f(x)g(x)$ is
$\sum_{k+\ell=n}{n\choose k}f_k g_\ell$. This is interpreted
as counting two (ordered) boxes worth of combinatorial objects:
\begin{itemize}
\item[(I)] the first of the type enumerated by $f(x)$, and 
\item[(II)] the second of $g$'s type;
\item[(III)] with the labels
distributed to the boxes in all possible ways. 
\end{itemize}
Note that $e^x-1$ is the series for \emph{unordered} collections
of languages. Thus, $(e^x-1)^k$ is the series for classifications
into $k$ \emph{ordered} families -- a $k!$ factor overcount of $S(n,k)$. 

Expanding the right hand side of (\ref{eqn:expgen}) by the binomial theorem, one deduces 
\begin{equation}
\label{eqn:stirling}
S(n,k)=\frac{1}{k!}\sum_{i=0}^k(-1)^i{k\choose i}(k-i)^n,
\end{equation}
where $S(n,0)=0$. One can lament that while this expression is explicit and nonrecursive, 
it is not manifestly nonnegative or even rational, even though it computes $S(n,k)\in {\mathbb Z}_{\geq 0}$.

Greenberg's problem for $n$ languages is computed by the {\bf Bell number} defined by
$B(n)=\sum_{k=0}^n S(n,k)$.
By (\ref{eqn:stirling}) we have a non-recursive expression for $B(n)$. However, computing
$B(25)$ by hand this way is impractical.

Actually, Greenberg's footnote on page 43 of \cite{Greenberg2} shows he knew the recurrence
\begin{equation}
\label{eqn:Bell}
B(n+1)=\sum_{k=0}^n {n\choose k}B(k),
\end{equation}
presumably used in his computation.\footnote{Greenberg cites \cite{Ore} that cites 
\cite{Epstein} (which has values of $B(n)$ for $n\leq 20$).
By 1962,
\cite{Levine} gave values for $n\leq 74$ and 
cite \cite{Gupta} who had values for 
$n\leq 50$. Hence $B(25)$ was 
known at least to experts, if not widely available, by the time of \cite{Greenberg2}. Anyway,
it seems likely he just computed it himself.}  
A combinatorial proof of (\ref{eqn:Bell}) follows by
decomposing all classifications of $n+1$ languages by the number of 
other languages in that family. Still, executing the calculation for $B(25)$ would have been 
a task.


Now, since $(e^x-1)^k/k!$ is the generating
series for $S(n,k)$, summing over all disjoint cases
$k$ of the number of families,
$\sum_{k=0}^{\infty} \frac{(e^x-1)^k}{k!}=e^{e^x-1}$
is the generating series for $B(n)$. 

Often, this is where the classroom or textbook analysis stops, with the computation treated as moot. However, carrying it out, using, e.g., {\tt Maple}, the instructor can make clear
the speed one calculates the answer. 
For example, in {\tt Maple} one computes as follows:

\begin{verbatim}
> 25!*coeftayl(exp(exp(x)-1),x^25);
                              4638590332229999353
\end{verbatim}
In the interpretation for the Joseph Greenberg problem, this is a small surprise.

\section{Back-of-the-envelope calculations}

An advantage of generating series is that they are adaptable to related enumerations. While Greenberg used $B(25)$ to support Mass comparison, similar numerics do not seem to be
consistent with a key consequence, his work on languages of the Americas \cite{Greenberg:americas}. 

\subsection{A simple plausibility test for Greenberg's classification}
Greenberg \cite{Greenberg:americas} classified $650$ indigenous languages of North and 
South America into three families\footnote{The $650$ figure is found, e.g., in 
\cite[pg.105]{keythink}}. This has been criticized for, e.g., not providing sufficient statistical evidence for claimed commonalities between languages. Specific to this situation is the use of a disputed method (Mass comparison) \emph{and} 
a large disagreement (two orders of magnitude) among linguists as to the number of families.
Side-stepping well-established lines of debate, imagine a restart using combinatorics. Roughly, how many language families should there be? 

Consider each classification in an unbiased way. 
The generating series for language families with between $a$ and $b$ families is
$\sum_{k=a}^b \frac{(e^x-1)^k}{k!}$.
Using this, the probability a random classification on $650$ languages has at most $3$ families is $0.238\times 10^{-843}\%$. (This is comparable to the probability of
randomly finding a prechosen atom from the observable universe correctly, ten times
in succession.) This hardly disproves Greenberg's classification -- but it does quantify how improbable it is from the baseline.

The probability that the number of families are in the ranges
$[50,110]$, $[111,120]$, $[121,130]$, $[131,140]$ and $[141,150]$ are 
$0.0000565\%$, $0.56\%$, $37.1\%$, $58.8\%$ and $3.5\%$ 
(rounded), respectively.
Hence almost all of the density is in the range $[121,150]$.
Thus, we would na\"ively guess (taking into account margin of error) that 
there are a few hundred families. This is decent agreement with today's consensus of $150$ to $180$ families \cite{Wiki, Wikib}.\footnote{For $1,000$ 
languages, which seems to me to be the upper end of
the number of described indigenous languages of the Americas, the model 
predicts the range $[180,200]$.}

\subsection{Language isolates}
A \emph{language isolate} is a language family consisting of only one language. Euskara, the ancestral language of the
Basque people,  is one such example. 

If all (modern) human languages originate from a single
source, then one should consider, as Greenberg does, 
language classifications with no isolates. Indeed, in his Americas classification,
Greenberg placed many generally regarded isolates (Yuchi, Chitimacha, Tunica, among others) into his Amerind superfamily.

Thus, suppose we consider such classifications for eight languages. 
The number of possibilities is the coefficient
of $\frac{x^8}{8!}$ in $e^{e^x-x-1}$ (which is $715$). While this is not small, it says that the
probability of a random language grouping having no language isolate is $\frac{715}{4140}$ or
approximately $17\%$. 
For $25$ languages, the probability is about $8.75\%$, whereas for $200$ languages it is
about $1.93\%$.  In the case of $n=500$ (about the number of
languages analyzed in \cite{Greenberg:africa}), it is $0.927\%$. For $650$ languages (roughly the number studied
in \cite{Greenberg:americas}) it is $0.747\%$. 
%
In other words, most random language family configurations have
a language isolate. This is not supportive of Greenberg's hypothesis.

By the product rule, the number of classifications on $n$ languages with $f$ families and $i$ isolates is the
coefficient of $\frac{x^{n}}{n!}$ in $\frac{x^i}{i!}\frac{(e^x-x-1)^{f-i}}{(f-i)!}$. For $n=650$, if $f=150$,
the expected number of isolates is about $9$ and nearly $5\%$ have more than $14$ isolates in the right tail. 
If instead $f=180$ (the upper range of the consensus number of families), the expected number of isolates is roughly $19.5$
with nearly $5\%$ having more than $26$ isolates in the right tail. 
Thus, the couple of dozen isolates/unclassified languages identified in the actual classification seems somewhat high (or our estimate seems somewhat low).
This would predict that a number of current isolates/unclassified languages should amalgamate, reducing the total number of families down towards $150$.
We refrain from a ``just-so'' argument for why the number of isolates seems higher than the model predicts.

\subsection{What about Africa?}
About $1,500$ languages (\cite[pg.104]{keythink}) of Africa were 
classified by Greenberg \cite{Greenberg:africa}, using Mass comparison,
into $4$ language families. Our na\"{i}ve tests 
also reject this conclusion, estimating instead a few hundred language families. 

	Yet Greenberg's classification is generally regarded as a success! 
However, in \cite[pg 561]{Sands} this success is qualified:
``for the majority of Africa's best documented languages''. In \emph{loc. cit.} it is noted that there is a documentation problem for less popular languages, and the total
number of languages in Africa varies substantially: from $2058$ by one count, to $1441$ in another.  Moreover, the review 
\cite{review} argues that there are $19$ language families, given present knowledge. 
Also, D.~Ringe suggested to us that Africa is special because of $2,000$ years of Bantu expansion that wiped out many language families. In any case, the situation is more complicated than first supposed.

Indeed, it is reasonable to tune the model, by factoring out uncontroversial language classifications, such as the 
Afro-Asiatic languages. However,
since an accounting of the current wisdom on African languages
seems to me to be a matter for an expert,
and would take us too far from the mathematics, further discussion may appear elsewhere. 

In conclusion, further interpretation, analysis and simulation may come 
from ideas in both combinatorics and comparative linguistics.

\section*{Acknowledgements}
We thank Donald Ringe for a helpful and encouraging email correspondence (and his appearance in the documentary), from which I learned much of the comparative linguistics 
background I've presented (any errors are my own). We also l
thank Eugene Lerman, Oliver Pechenik (especially for suggesting a criticism of \cite{Greenberg:americas}), Anh Yong, and the Fall 2013 Math 580 class at UIUC.
The author was supported by an  NSF grant.

\end{document}